# Implicit multiderivative collocation solvers for linear partial differential equations with discontinuous Galerkin spatial discretizations

Jochen Schütz[*1], David C. Seal[†2], and Alexander Jaust[‡1]

[1]*Faculty of Sciences, Hasselt University, Agoralaan Gebouw D, BE-3590 Diepenbeek*
[2]*United States Naval Academy , Department of Mathematics , 572C Holloway Road, Annapolis, MD 21402, USA*

In this work, we construct novel discretizations for the unsteady convection-diffusion equation. Our discretization relies on multiderivative time integrators together with a novel discretization that reduces the total number of unkowns for the solver. These type of temporal discretizations come from an umbrella class of methods that include Lax-Wendroff (Taylor) as well as Runge-Kutta methods as special cases. We include two-point collocation methods with multiple time derivatives as well as a sixth-order fully implicit collocation method that only requires a total of three stages. Numerical results for a number of sample linear problems indicate the expected order of accuracy and indicate we can take arbitrarily large time steps.

## 1 Introduction

We consider discretizations of the scalar convection-diffusion equation on a two-dimensional domain $\Omega \subset \mathbb{R}^2$:

$$\partial_t w + \nabla \cdot (\mathbf{c} w - \varepsilon \nabla w) = g, \qquad \forall (x, t) \in \Omega \times (0, T_{end}], \tag{1}$$

$$w(x, 0) = w_0(x) \qquad \forall x \in \Omega, \tag{2}$$

---

[*]Electronic address: `jochen.schuetz@uhasselt.be`
[†]Electronic address: `seal@usna.edu`
[‡]Electronic address: `alexander.jaust@uhasselt.be`

where $\mathbf{c} \in \mathbb{R}^2$, $\varepsilon \in \mathbb{R}$ and $g \in L^2(\mathbb{R}^2 \times \mathbb{R}^+)$ are prescribed parameters. Our aim is to apply a classical discontinuous Galerkin (DG) method to the spatial part of (1), and to advance the solution in time with an implicit multiderivative Runge-Kutta method [25]. The novelty in this work is the combination of these implicit multiderivative methods with the discontinuous Galerkin method, and the result is that we are able to take arbitrarily large time steps with high-order solvers all the while reducing the total number of stages that would normally be required to obtain the same order of accuracy.

## 1.1 Background

This work makes use of multistage multiderivative (MSMD) time integrators, which are best described in the context of ordinary differential equations (ODEs). For an ODE defined by $\partial_t y = f(y)$, multiderivative schemes make use of higher temporal derivatives of the unknown solution $y$. These time derivatives can be expressed recursively in terms of $f$ and its derivatives. For example, the second time derivative can be expressed as $\partial_{tt} y = f'(y)\partial_t y = f'(y)f(y)$, where $f'(y)$ is the Jacobian of $f$ with respect to $y$. These higher derivatives, together with additional stages, form part of the foundation of all MSMD methods, we well as their even more general multistep-multistage-multiderivative extensions [23].

Because additional information is fed into the algorithm, MSMD methods can be constructed to obtain higher order of accuracy than a standard Runge-Kutta scheme with the same number of stages. For example, with a total of $s$ stages, Butcher showed that it is possible to construct a Runge-Kutta method that obtains $(2s)^{th}$-order accuracy [9]. On the other hand, Stroud and Stancu [47] have shown that it is possible to obtain a method with a total of $(2(s+M))^{th}$-order accuracy, when a total of $M$ (even) derivatives of the right hand side function are considered. This is because the additional degrees of freedom required to obtain higher order accuracy can be found by searching for higher derivatives in place of adding additional stages. Methods from this class can be constructed as explicit or implicit solvers, and the implicit solvers can be designed in a way such that they fulfill desirable properties, such as $A$-stability [16]. A special case of these solvers include all Taylor methods, where a total of one stage is considered, and the coefficients of the higher derivatives are picked from the Taylor series of the solution. These can also be used to construct implicit or explicit solvers. Note that introducing a time history of the solution can serve as an alternative to increasing the order of a (multistage) Runge-Kutta method. This defines the class of so-called general linear methods [10, 8, 50], that can be thought of as a multistage multistep method.

In the context of PDEs, the special case of a Taylor discretization is typically called a Lax-Wendroff solver. This is attributed to the original work of Lax and Wendroff from 1960 where they construct a second-order solver by incorporating the second derivative of the PDE into their method [33]. More recently, higher order (i.e., solvers with order greater than two) versions of these solvers have been investigated for finite volume [26], finite difference [41, 32, 46, 12, 14, 52], and discontinuous Galerkin discretizations [40, 22, 36]. A large community centered around Arbitrary DERivative (ADER) discretizations has been very successful with constructing arbitrary order explicit solvers for hyperbolic problems in this category [48, 51, 15, 7], and much of their work relies on symbolic software to generate their code base.

Although seldom used, the multistage multiderivative methods have been investigated (for ODEs) since as early as the 1960's for problems in celestial mechanics [43], and later on for various other differential equations [20, 19]. The multistage multiderivative flavor of these solvers have only recently attracted attention as a mechanism for discretizing partial differential equations [49, 45]. In [13] it is shown that the multistage multiderivative formulations can be constructed to contain the so-called strong stability preserving property, and these solvers are currently being investigated as useful time discretizations for equations of gas dynamics [34, 39, 38]. An extensive review of these methods can be found in [45].

In our previous work [30], we began an attempt to couple DG and two-point two-derivative methods. This earlier work was based on the Cauchy-Kovalevskaya procedure [42], which means that one takes the original PDE, and expresses the temporal derivatives of the unknown $w$ in terms of the spatial derivatives. Based on this ansatz, we introduced additional variables modeling the spatial derivatives of $w$ to model $w_{tt}$. This term is used to express the spatial derivatives of $w$ up to order four. This procedure has lead to a method that, although having quite large linear systems of equations to solve, is in runtime comparable to classical time integration schemes. Note that because the scheme was based on the Cauchy-Kovalevskaya procedure together with similar tricks used to define Lax-Wendroff discontinuous Galerkin solvers, our previous scheme is *not identical* to applying a two-derivative scheme to the ordinary differential equation that results from the method-of-lines formulation stemming from the DG spatial discretization of (1). This lead to unwanted features such as a sometimes quite severe loss of stability reducing an implicit solver to finite time steps on the same order as an explicit method. The present work is directed at mending this undesirable feature.

## 1.2 Summary of work

In this work, we redirect our efforts away from discretizing each higher derivative term separately, and instead construct a discretization that ends up being equivalent to the method-of-lines (MOL) formulation of the partial differential equation (PDE). That is, in place of attempting to do anything special to define higher derivatives, which is commonplace with Taylor type discretizations, we construct a solver that is equivalent to applying the multiderivative methods directly to the MOL discretization of the PDE. That is to say, instead of following the common practice of defining higher derivatives by differentiating basis functions (or using reduced order stencils), we instead take great care to construct these higher derivatives in such a way that they end up being identical to differentiating the large system of ODEs defined by the MOL discretization of the PDE. We make this important modification in order to construct and prove the stability of our solver. In doing so, one obvious complaint is that this has the potential to increase the size of the effective stencil of the method. We address this issue by introducing only *one* additional variable that is used to define all of the mixed derivatives of the solution.

In this work, we consider a total of two different types of implicit collocation methods,

and we couple each of these with the classical discontinuous Galerkin discretization of the PDE. The methods we consider can be classified into a total of two separate categories:

- **Two-point multiderivative collocation methods**. These methods use a total of two quadrature points (one at a known time value $t = t^n$, and another at the next time value $t = t^{n+1}$). They reach high order accuracy by increasing the number of derivatives (we consider methods with a total of three derivatives of the unknown) that are evaluated at each time point.

- **Fully implicit multiderivative collocation methods**. These methods increase the order of accuracy by increasing the number of stages, and they too can increase order by including higher derivatives. These methods can be constructed by first defining a set of collocation (quadrature) points, fitting a Hermite-Birkhoff interpolant, and then integrating the result. Because a two-derivative method with three quadrature points can obtain a total of sixth-order accuracy, we only include a single method from this category. Because sixth-order accuracy is very high order, in this work we do not pursue adding additional stages and point to this solver as a proof of concept.

Both of these classes of methods fall under the umbrella category of multiderivative Runge-Kutta methods. This broader category also encompasses all classical Runge-Kutta solvers, as well as Taylor (or Lax-Wendroff) methods, but not all of them are $A$-stable. Here, we only consider $A$-stable solvers, which we show is an important property that we leverage to define stable numerical discretizations for the PDE.

The paper is structured as follows: in Section 2, we provide a brief review of the discontinuous Galerkin discretization. This section serves to the notation that is used throughout the remainder of this work. In Section 3, we review classical multiderivative time discretizations in the context of ordinary differential equations, and in Section 4, we couple the two discretizations to define the new fully discrete solver. Finally, in Section 5 we present numerical results, and in Section 6 we wrap up with conclusions and point to future work.

## 2 Discontinuous Galerkin discretizations in space

We begin by introducing the DG discretization that we use in this sequel. In principle, one can substitute the scheme by one's favorite scheme, as long as it is coercive (and therefore stable) and in primal form.

Before introducing the scheme in detail, we shortly define the (rather standard) notation needed. Based on a triangulation of the domain $\Omega \subset \mathbb{R}^2$ into

$$\Omega = \bigcup_{k=1}^{N_e} \Omega_k, \tag{3}$$

we define the standard space of broken polynomials,

$$V_h := \{f \in L^2(\Omega) \mid f_{|\Omega_k} \in \Pi^p(\Omega_k)\},$$

where the set $\Pi^p(\Omega_k)$ is the space of all polynomials on $\Omega_k$ having total degree $p$. Other choices are possible, in particular, an adaptive polynomial degree does not pose any particular problems. We do not pursue this any further in this work.

Cell-wise integration over $\Omega$ is denoted by the scalar product $(\cdot,\cdot)$, while $\langle\cdot,\cdot\rangle$ denotes edge-wise integration over the skeleton of the triangulation. To consider functions on an edge $e_{k,l} := \Omega_k \cap \Omega_l$, $k \neq l$, we need to define 'inner' and 'outer' values. Let $e_{k,l}$ be equipped with a normal vector $\mathbf{n}$, and let $x \in e_{k,l}$, then we define for a function $\varphi_h \in V_h$,

$$\varphi_h^\pm := \lim_{\delta \to 0} \varphi_h(x \pm \delta\mathbf{n}).$$

Furthermore, we define average and jump, respectively, as

$$\{\varphi_h\} := \frac{\varphi_h^- + \varphi_h^+}{2}\mathbf{n}, \qquad [\![\varphi_h]\!] := \left(\varphi_h^- - \varphi_h^+\right)\mathbf{n}.$$

Note that, while $\varphi_h^\pm$ depends on the orientation of the normal, jumps and average do not.

To discretize the convective terms, we employ a standard upwinding technique. To this end, we define

$$w_h^{up} := \begin{cases} w_h^-, & \mathbf{c} \cdot \mathbf{n} > 0, \\ w_h^+, & \text{otherwise.} \end{cases}$$

Note that this definition is independent of the orientation of $\mathbf{n}$. We discretize the viscous term with a symmetric interior penalty method [3].

With these definitions in place, the discretization of (1) in space yields the task of seeking $w_h \in C^0([0, T_{end}], V_h)$, such that

$$(\partial_t w_h, \varphi_h) - (\mathbf{R}(w_h), \nabla\varphi_h) + \langle \mathbf{R}_e(w_h; \varphi_h)\rangle = (g, \varphi_h), \quad \forall \varphi_h \in V_h, \tag{4}$$

where $\mathbf{R}(w_h) \equiv \mathbf{R}(w_h, \nabla w_h)$ and $\mathbf{R}_e(w_h; \varphi_h) \equiv \mathbf{R}_e(w_h, \nabla w_h; \mathbf{n}; \varphi_h, \nabla\varphi_h)$ are the cell and edge discretizations of DG-type given by

$$\mathbf{R}(w_h) := \mathbf{c}w_h - \varepsilon\nabla w_h, \tag{5}$$

$$\mathbf{R}_e(w_h; \varphi_h) := -\{\mathbf{c}w_h^{up}\} \cdot [\![\varphi_h]\!] + \varepsilon\left(\{\nabla w_h\} - \frac{\eta}{h_e}[\![w_h]\!]\right) \cdot [\![\varphi_h]\!]$$
$$+ \varepsilon[\![w_h]\!]\{\nabla\varphi_h\}. \tag{6}$$

Note that these operators are linear in $w_h$ and $\varphi_h$. The value $h_e$ is the length of an edge, and $\eta$ is a user-defined parameter that must be positive and larger than a certain threshold, (see [3] for details). Upon inserting a basis for $V_h$ into (4), we rewrite this as a large complex linear system of differential equations,

$$\partial_t \mathbf{w}_{DG} = \mathbf{A}_{DG}\mathbf{w}_{DG} + \mathbf{b}_{DG}, \tag{7}$$

where $\mathbf{w}_{DG}$ is a vector of unknowns, $\mathbf{A}_{DG}$ is a matrix representing the difference operators, and $\mathbf{b}_{DG}$ is a vector representing the source terms.

**Remark 1.** *Note that for $\varepsilon > 0$ and polynomial order $p = 0$, the method presented above is not meaningful. For the case where diffusion appears, we will therefore not show any results produced with piecewise constant ansatz functions. However, the methodology we present does not rely on the particular choice of the DG discretization, it is therefore possible to substitute $\mathbf{R}$ and $\mathbf{R}_e$ by other stable discretization types which can handle the $p = 0$ case for diffusion.*

The proof of stability of our solver relies on the following lemma (whose proof can be found in [3, 28, 27]):

**Lemma 1.** *Let the triangulation be conforming and shape-regular; and let $p > 0$. Then, there exist an $\eta^* > 0$, such that for all $\eta > \eta^*$ and $\varphi_h \in V_h$, we have*

$$-(\mathbf{R}(\varphi_h), \nabla \varphi_h) + \langle \mathbf{R}_e(\varphi_h; \varphi_h) \rangle \geq 0. \tag{8}$$

The above statement means that under the conditions mentioned in the Lemma, the matrix $\mathbf{A}_{DG}$ is negative semi-definite (for $\varepsilon > 0$, it is negative definite). This means that the real part of each eigenvalues of the matrix is negative:

**Corollary 1.** *Let $\lambda$ be an eigenvalue of $\mathbf{A}_{DG}$ under the conditions of Lemma 1. Then, the real part of lambda is negative. That is, if $\lambda$ is an eigenvalue of $\mathbf{A}_{DG}$, then $\Re(\lambda) \leq 0$. Furthermore, if $\varepsilon > 0$, then $\Re(\lambda) < 0$.*

When $\mathbf{c} \neq 0$, the matrix is not symmetric, and therefore the statement of this Corollary is not directly clear. Here, we prove the case where $\varepsilon > 0$. The case where $\varepsilon = 0$ is nominally different.

*Proof.* Thanks to Eqn. (8), we have that $\mathbf{m}^T \mathbf{A}_{DG} \mathbf{m} < 0$ for all non-zero vectors $\mathbf{m}$. Furthermore, because $\mathbf{A}_{DG} \mathbf{v} = \lambda \mathbf{v}$, we have $\mathbf{A}_{DG} \overline{\mathbf{v}} = \overline{\lambda} \overline{\mathbf{v}}$ because $\mathbf{A}_{DG}$ has real coefficients. Therefore, we can write

$$\mathbf{A}_{DG} \Re(\mathbf{v}) = \Re(\lambda) \Re(\mathbf{v}) + \Im(\lambda) \Im(\mathbf{v}),$$
$$\mathbf{A}_{DG} \Im(\mathbf{v}) = -\Im(\lambda) \Re(\mathbf{v}) + \Re(\lambda) \Im(\mathbf{v}).$$

(Note that $\Re$ and $\Im$ denote real and imaginary part, respectively.) We multiply the first line from the left with $\Re(v)^T$, and the second line with $\Im(v)^T$, and further exploit the negative semi-definiteness of the operators to observe

$$0 > \Re(\lambda)\|\Re(\mathbf{v})\|^2 + \Im(\lambda)\Re(\mathbf{v})^T \Im(\mathbf{v}),$$
$$0 > -\Im(\lambda)\Im(\mathbf{v})^T \Re(\mathbf{v}) + \Re(\lambda)\|\Im(\mathbf{v})\|^2.$$

Upon adding terms, we observe

$$0 > \Re(\lambda) \left( \|\Re(\mathbf{v})\|^2 + \|\Im(\mathbf{v})\|^2 \right),$$

which proves the claim. □

**Remark 2.** *Corollary 1 is the reason we favor A-stable schemes.*

With these preliminaries out of the way, we are now prepared to discuss the various temporal discretizations that we use in this work. We begin with a description of how these solvers operate on ordinary differential equations.

## 3 Multiderivative discretization in time

In this section, we briefly review multiderivative Runge-Kutta methods. In the subsequent section, we apply the temporal discretization of (4) to these solvers.

An $M$-derivative Runge-Kutta solver is defined by a total of $M$ Butcher tableaux $\{a^{(1)}, a^{(2)}, \ldots, a^{(M)}\}$, each of size $s \times s$, where $s$ refers to the total number of stages of the solver, and $M$ is the total number of derivatives under consideration. The internal stages of an $M$-derivative Runge-Kutta scheme are defined as

$$y_{(i)}^n = y^n + \sum_{m=1}^{M} \Delta t^m \sum_{j=1}^{s} a_{ij}^{(m)} \partial_t^m y_{(j)}^n, \quad i = 1, \ldots, s, \tag{9a}$$

and the final update is given by

$$y^{n+1} = y^n + \sum_{m=1}^{M} \Delta t^m \sum_{i=1}^{s} b_i^{(m)} \partial_t^m y_{(i)}^n. \tag{9b}$$

The two-point collocation schemes we consider have a total of $s = 2$ stages, and up to $M = 3$ time derivatives of the unknown. The fully implicit collocation method we consider has a total of $s = 3$ stages, and $M = 2$ derivatives of the unknown.

### 3.1 Two-point multiderivative collocation methods

The first class of solvers we consider are two-point multiderivative collocation methods. These methods have a total of $s = 2$ stages, and the abscissa are sampled at times $t = t^n$ and $t = t^{n+1}$ only. These methods can be derived from the rational that one 'prescribes' $k$ derivatives at time $t^n$ and $l$ derivatives at time $t^{n+1}$. That is, we first fit Hermite-Birkhoff interpolant [29, 37] to the unknown function $y(t)$, and then integrate the result to define the solver. This results in the following explicit expression for the numerical solver:

$$\sum_{j=0}^{m} \Delta t^j (\partial_t^j y)^{n+1} P^{(m-j)}(0) = \sum_{j=0}^{m} \Delta t^j (\partial_t^j y)^n P^{(m-j)}(1), \tag{10}$$

where $P(t) = \frac{t^k(t-1)^l}{(k+l)!}$ and $\partial_t^j y$ is the $j$−th temporal derivative of the solution $y$ to the ODE [24]. In this work, we employ two-point schemes for two and three derivatives, respectively. Each of these schemes can be written in the form

$$\begin{aligned}\frac{y^{n+1} - y^n}{\Delta t} =\ & \left(\alpha_1 \partial_t y^n - \beta_1 \partial_t y^{n+1}\right) + \Delta t \left(\alpha_2 \partial_t^2 y^n - \beta_2 \partial_t^2 y^{n+1}\right) \\ & + \Delta t^2 \left(\alpha_3 \partial_t^3 y^n - \beta_3 \partial_t^3 y^{n+1}\right), \end{aligned} \tag{11}$$

for some values of $\boldsymbol{\alpha} = (\alpha_1, \alpha_2, \ldots)$ and $\boldsymbol{\beta} = (\beta_1, \beta_2, \ldots)$. The values that we use in this work are summarized in Table 1, where we categorize the solvers based on the order of the method.

Table 1: Two-point multiderivative schemes. Those with $\alpha_3 = \beta_3 = 0$ only need two derivatives, which means additional variables do not need to be computed.

| Order | $(k,l)$ | $\alpha_1$ | $\alpha_2$ | $\alpha_3$ | $\beta_1$ | $\beta_2$ | $\beta_3$ |
|---|---|---|---|---|---|---|---|
| 3 | $(1,2)$ | $1/3$ | $0$ | $0$ | $-2/3$ | $1/6$ | $0$ |
| 4 | $(2,2)$ | $1/2$ | $1/12$ | $0$ | $-1/2$ | $1/12$ | $0$ |
| 5 | $(2,3)$ | $2/5$ | $1/20$ | $0$ | $-3/5$ | $3/20$ | $-1/60$ |
| 6 | $(3,3)$ | $1/2$ | $1/10$ | $1/120$ | $-1/2$ | $1/10$ | $-1/120$ |

**Remark 3.** *The values in Table 1 are very much related to Padé approximations for the exponential function. As an example, the Padé-approximation of the exponential function of order $(2,3)$ is given by*

$$P_{2,3}(z) = \frac{1 + 2/5z + 1/20z^2}{1 - 3/5z + 3/20z^2 - 1/60z^3}.$$

*These coefficients are identical to those found in the fifth-order scheme defined in Table 1.*

The following very important Lemma comments on the stability of these schemes, and is critical for us to define a stable numerical solver:

**Lemma 2.** *All of the methods shown in Table 1 are A-stable.*

*Proof.* Thanks to Remark 3, all of these methods are on the first or second subdiagonal in the set of Padé approximations to $e^z$. In [17], it is shown that each these entries are $A$-stable. Moreover, the methods on the sub-diagonal exhibit stiff decay, and therefore those methods are $L$-stable. □

### 3.2 Fully implicit multiderivative collocation methods

Extensions of two-point collocation methods can proceed in several directions. For example, it is possible to add more stages, more derivatives, or more steps to the solver to increase the order of accuracy [24]. For example, the so-called multistep multiderivative methods [18, 21, 31] increase the order of accuracy by adding a time history of the solution. In this work, we prefer to restrict our attention to self-starting single-step methods, and therefore we only consider additional stages in order to improve the order of accuracy of the method. Again, a fully implicit multiderivative collocation method can be easily derived by first fitting a Hermite-Birkhoff interpolant to the unknown function $y(t)$, and then integrating the result. The resulting scheme takes the form of a multiderivative Runge-Kutta method.

One nice property is that by deriving fully implicit multiderivative solvers in this manner we automatically know that they satisfy the correct order conditions. This is a result of Obreshkov's formula [35], which can be thought of as a generalization of Taylor's theorem to Hermite-Birkhoff interpolation. For general Runge-Kutta methods, this is normally a non-trivial task to accomplish, and more to the point, finding $A$-stable methods is already a difficult enough task to accomplish, let alone trying to find one that is high-order.

As a proof of concept, we restrict our attention to a single method that we derive in this manner. More specifically, the method we consider is a fully implicit two-derivative collocation method that we find to be $A$-stable. The development of higher order $A$-stable multistage multiderivative methods is reserved for future work.

Based on the work in [47], we compute the coefficients for a sixth-order scheme involving $s = 3$ stages. The collocation points are at the (normalized) time instances $\mathbf{t} = (0, 1/2, 1)^T$, which is important to know when considering non-autonomous differential equations. The Butcher tableau $a^{(1)}$ (for first derivative) and $a^{(2)}$ (for second derivative) read as follows:

$$a^{(1)} = \begin{pmatrix} 0 & 0 & 0 \\ 101/480 & 8/30 & 55/2400 \\ 7/30 & 16/30 & 7/30 \end{pmatrix}, \quad a^{(2)} = \begin{pmatrix} 0 & 0 & 0 \\ 65/4800 & -25/600 & -25/8000 \\ 5/300 & 0 & -5/300 \end{pmatrix}. \tag{12}$$

Because there are a total of six pieces of information used to define the Hermite-Birkhoff fit for the right hand side function (two each at times $t^n$, $t^{n+1/2}$, and $t^{n+1}$), this method achieves a total of sixth-order accuracy. More specifically, the error in this approximation can be found from Obreshkov's formula [35]. We stop to point out that one advantage of this solver is that it does not require extra derivatives of the unknown function, but this comes at the expense of adding an additional stage.

**Lemma 3.** *The method defined by Butcher tableaux (12) is $A$-stable.*

The proof of this Lemma is straightforward and is omitted for brevity.

## 4 The fully discrete solver

With these preliminaries out of the way, we are now prepared to describe the fully discrete solver proposed in this work. In a straightforward way, one can of course differentiate (7) and obtain an explicit expression for $\partial_{tt}\mathbf{w}_{DG}$. Together with Lemma 1 and the A-stability of the methods involved, this yields a stable algorithm.

**Lemma 4.** *With the A-stability of the involved time discretization schemes, and the coercivity of $\mathbf{A}_{DG}$, see Lemma 1 and Corollary 1, the application of those time integration schemes to the DG semi-discretization (7) yields a stable scheme.*

As the matrix $\mathbf{A}_{DG}^2$ occurs in the explicit representation of $\partial_{tt}\mathbf{w}_{DG}$, this involves an increase of the stencil of the method. To keep the compact stencil of the DG method - which is one of its many advantages - we introduce the additional variable $\sigma_h$ (which is assumed to have the same dimensionality as $w_h$, thus is a scalar) that fulfills

$$(\sigma_h, \varphi_h) - (\mathbf{R}(w_h), \nabla\varphi_h) + \langle \mathbf{R}_e(w_h; \varphi_h) \rangle = (g, \varphi_h) \quad \forall \varphi_h \in V_h. \tag{13}$$

Note that this definition is very similar to (4), in fact, only the first term differs.

With the help of this variable, we can express $\partial_{tt}w_h$ as follows:

**Lemma 5.** *Let $\sigma_h$ be defined as in (13). Then, for $w_h$ as defined in (4), there holds*

$$(\partial_{tt} w_h, \varphi_h) = (\mathbf{R}(\sigma_h), \nabla \varphi_h) - \langle \mathbf{R}_e(\sigma_h; \varphi_h) \rangle + (\partial_t g, \varphi_h). \tag{14}$$

*Proof.* As both $\mathbf{R}$ and $\mathbf{R}_e$ are linear functions, the vector $\mathbf{w}_{DG}$ containing the basis coefficients of $w_h$ for a certain set of basis functions fulfills

$$\partial_t \mathbf{w}_{DG} = \mathbf{A}_{DG} \mathbf{w}_{DG} + \mathbf{b}_{DG}, \tag{7}$$

and its second derivative can be easily computed as

$$\partial_{tt} \mathbf{w}_{DG} = \mathbf{A}_{DG} \left( \mathbf{A}_{DG} \mathbf{w}_{DG} + \mathbf{b}_{DG} \right) + \partial_t \mathbf{b}_{DG}.$$

By the construction of $\sigma_h$, we know that its associated vector $\boldsymbol{\sigma}_{DG}$ of basis coefficients fulfills $\boldsymbol{\sigma}_{DG} = \mathbf{A}_{DG} \mathbf{w}_{DG} + \mathbf{b}_{DG}$. Applying $\mathbf{A}_{DG}$ once more, done in (14) on the right-hand side, completes the proof. $\square$

**Remark 4.** *Note that, as said earlier, $\mathbf{A}_{DG}^2$ is an operator that has a larger stencil than $\mathbf{A}_{DG}$ has. In fact, one has to incorporate neighbor-neighbors. If one defines $\sigma_h \approx \partial_t w_h$ as auxiliary variable instead and expresses $\partial_{tt} w_h$ as in (14), the stencil is not enlarged. Furthermore, the assembly process is quite simple, as the terms used for $\sigma_h$ are the same as those for $w_h$, so most computations must be done only once.*

**Remark 5.** *Note that the proof of Lemma 5 already gives a glimpse on how to obtain a scheme in the case of a nonlinear operator $\mathbf{A}_{DG}(\mathbf{w}_{DG})$. For example, when discretizing the Navier-Stokes equations: Use the same trick on $\sigma_h$, and then multiply it with the derivative of $\mathbf{A}_{DG}$, evaluated at $\mathbf{w}_{DG}$. This will include some tedious differentiation operations that, however, can be simplified by using automatic differentiation. Also, tricks similar to those recently performed in [52], where higher derivatives of the flux function are approximated using high-order finite differences, can most certainly be applied here. The performance, analysis, and implementation of such a type of algorithm will be left for future work.*

**Remark 6.** *Of course with the same reasoning, one can introduce higher derivatives. In this work, we consider algorithms with up to three (temporal) derivatives. This means that beyond $\sigma_h$, there must be a third variable which we call $\tau_h$. It is defined in close analogy to (13) as*

$$(\tau_h, \varphi_h) - (\mathbf{R}(\sigma_h), \nabla \varphi_h) + \langle \mathbf{R}_e(\sigma_h; \varphi_h) \rangle = (\partial_t g, \varphi_h) \quad \forall \varphi_h \in V_h. \tag{15}$$

*Lemma 5 holds with obvious modifications.*

With these remarks, the algorithm is completely defined. As an example, the method for a two-point scheme can be written as follows:

$$\begin{aligned}\frac{\mathbf{w}_{DG}^{n+1} - \mathbf{w}_{DG}^n}{\Delta t} = & \mathbf{A}_{DG} \left( \alpha_1 \mathbf{w}_{DG}^n - \beta_1 \mathbf{w}_{DG}^{n+1} \right) + \Delta t \mathbf{A}_{DG} \left( \alpha_2 \boldsymbol{\sigma}_{DG}^n - \beta_2 \boldsymbol{\sigma}_{DG}^{n+1} \right) \\ & + \Delta t^2 \mathbf{A}_{DG} \left( \alpha_3 \boldsymbol{\tau}_{DG}^n - \beta_3 \boldsymbol{\tau}_{DG}^{n+1} \right),\end{aligned} \tag{16}$$

where we assume $\mathbf{b}_{DG} = \mathbf{0}$ for the sake of readability. Here, the definition of $\boldsymbol{\tau}_{DG}$ is similar to the ones of $\mathbf{w}_{DG}$ and $\boldsymbol{\sigma}_{DG}$, it denotes the vector with basis coefficients associated to $\tau_h$. Integrating the source term is straightforward, following the lines of (11).

# 5 Numerical results

In this section, we present numerical results, demonstrating the performance of the scheme. In all the results, we work with periodic boundary conditions to alleviate influence from the boundary. We note, though, that the algorithm can be easily adapted to handle non-periodic boundary conditions. In fact, the treatment of the boundary conditions is usually hidden in $\mathbf{A}_{DG}$ and $\mathbf{b}_{DG}$, and so the multiderivative algorithm does not need to explicitly deal with it. This is in contrast with Lax-Wendroff schemes, where of course the boundary conditions influence the representation of $\partial_{tt} w$.

The meshing tool Netgen [44] serves as a basis for the code. The linear systems of equations are solved through PETSc [6, 5, 4], either using a GMRES scheme that is converged up to a relative tolerance of $10^{-10}$ with ILU(2) preconditioner; or through a direct solver. The error $e_h$ is always defined as $L^2$−error at time $T_{end}$, i.e.,

$$e_h := \|w(t = T_{end}) - w_h(t = T_{end})\|_{L^2}.$$

## 5.1 Convection equation

We start with the relatively simple convection test case, characterized by a constant advection speed $\mathbf{c} = (1, 1)^T$, a zero diffusion coefficient $\varepsilon = 0$ and the initial conditions $w_0(x, y) = \sin(2\pi x) \sin(2\pi y)$. The domain $\Omega = [0, 1]^2$ is the unit cube. The solution is a transport of the initial conditions in direction $\mathbf{c}$, and we choose a final time of $T_{end} = 1.0$.

### 5.1.1 The convection equation: Two-point schemes

We begin with the third and fourth order two-point schemes defined by (11) and Tbl. 1. These schemes only need one additional temporal derivative, thus $\tau_h$ is not computed.

Convergence of the error $e_h$ versus time step size $\Delta t$ is plotted in Fig. 1 for various polynomials orders; the smallest mesh consists of two triangular elements, both spatial and temporal refinement is uniform. On the left-hand side, the third-order method is plotted, the right-hand side shows the fourth-order method. It can be seen that optimal convergence orders $\min\{p+1, 3\}$ and $\min\{p+1, 4\}$, respectively, are reached. We have performed this exercise for different values of the CFL number (smaller and larger than one), and the error plots look quite similar. It is only that for large CFL numbers that it takes longer until optimal order is reached due to the under-resolution of time. This, however, is common to all time integration schemes.

In contrast to some of the multiderivative test cases reported on in [30], we do not observe any stability problems. What seems to be quite relevant is the choice of the preconditioner. We find that using too simple of a preconditioner (e.g., Jacobi) results in a non-convergence of the linear system for CFL numbers above one. We find that using the ILU(2) choice is highly robust for the convection equation.

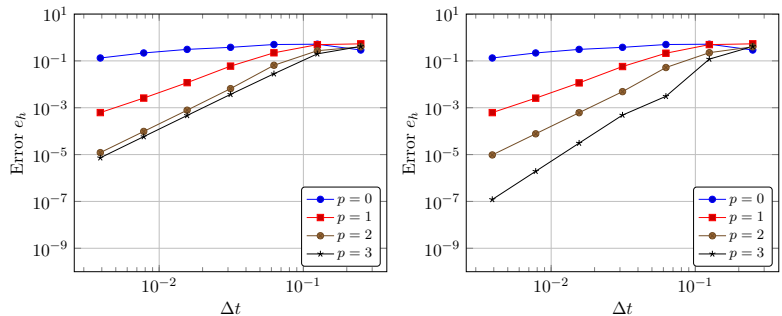

Figure 1: Numerical results for the convection equation with parameters $\mathbf{c} = (1,1)^T$, $T_{end} = 1.0$ and $\varepsilon = 0$. Two temporal derivatives of the DG scheme were needed. Time step size for the smallest triangular mesh, consisting of two elements, was chosen to be $\Delta t = 0.25$. Left plot: results using the third-order two-point scheme, see Tbl. 1; right plot: results using the fourth-order two-point scheme, see Tbl. 1.

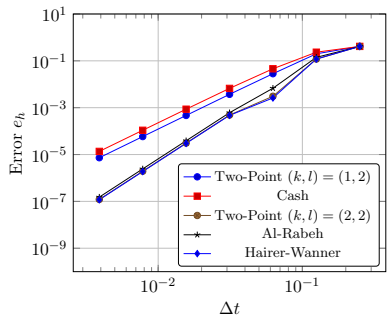

Figure 2: Numerical results for the convection equation with parameters $\mathbf{c} = (1,1)^T$, $T_{end} = 1.0$ and $\varepsilon = 0$. The polynomial order was chosen to be $p = 3$, and we compare multiple time integration schemes against each other. $\Delta t = 0.25$ was chosen on the coarsest mesh, consisting of two triangular elements.

In Fig. 2, we compare the performance of the third- and fourth-order two-point schemes to those of the more established schemes of diagonally implicit Runge-Kutta (DIRK) type. More precisely, we test the schemes for polynomial order $p = 3$ against the more or less classical DIRK schemes by Cash [11] (this scheme is also due to Alexander [2]), Al-Rabeh [1] and Hairer and Wanner [25]. The last two DIRK schemes are of order four, the first one is of order three. They consist of three, four and five stages, respectively. (Note that a two-point scheme formally consists of two stages, but the schemes we investigate also have additional derivatives.) The parameters are the same as those in Fig. 1. It is obvious that the two-point schemes behave as good (or sometimes even slightly better) than the classical DIRK schemes. With regard to stability, all methods perform equally. Computing equations with $\Delta t \gg \Delta x$ is done in a stable way, and also the error curves behave quite similarly.

### 5.1.2 The convection equation: Higher-order derivatives

We continue with the use of the three-derivative two-point schemes presented in Tbl. 1, being of order five and six, respectively. The reason we are studying these schemes is to test how well the method can be extended to higher orders, not only using stages but also using additional derivatives. Fig. 3 shows numerical results for the same test case presented above. As the methods are of higher order than before, we include polynomial orders up to $p = 5$. Again, optimal orders of convergence (here, it is $\min\{p+1, 5\}$ and $\min\{p+1, 6\}$, respectively) are obtained.

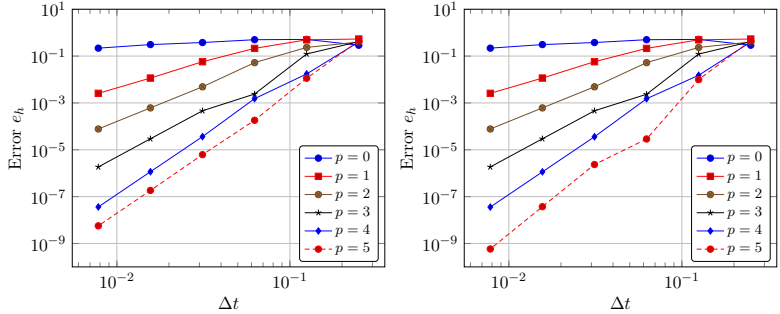

Figure 3: Numerical results for the convection equation with parameters $\mathbf{c} = (1, 1)^T$, $T_{end} = 1.0$ and $\varepsilon = 0$. Three temporal derivatives of the DG scheme are needed. Time step size for the smallest triangular mesh, consisting of two elements, was chosen to be $\Delta t = 0.25$. Left plot: results using the fifth-order two-point scheme, see Tbl. 1; right plot: results using the sixth-order two-point scheme, see Tbl. 1.

### 5.1.3 The convection equation: A three point scheme

We next test the utility of adding additional stages to the solver. Because a two-derivative method with a single additional stage can obtain sixth-order accuracy, we only consider one method of this type in this work. Ultimately, tests with the multiderivative collocation scheme given in (9) (with Butcher tableaux in (12)) are performed using the same parameters as above. Results can be seen in Fig. 4; the optimal order of $\min\{p+1, 6\}$ is achieved in all cases. Again no stability problems are observed.

We find that this method is competitive against other well-known implicit methods, such as the Gauß-Legendre or Radau methods. For example, in Fig. 5, we compare the sixth-order Gauß-Legendre Runge-Kutta method against the multiderivative collocation method used in this work. The polynomial order used is $p = 5$, so that one can indeed observe sixth order convergence for all the methods. With the same parameters that we use in the previous figures, we are unable to discern any difference between the schemes, which tells us that the spatial error is dominating. In order to elucidate the difference between the solvers we re-compute the test case with an initial time step of size $\Delta t = 1.0$ on the coarsest mesh and then refine from there. We observe that the multiderivative scheme shows a slight advantage in terms of the error.

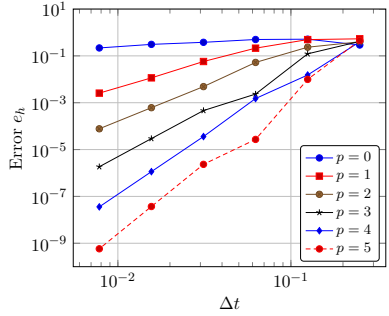

Figure 4: Numerical results for the convection equation with parameters $\mathbf{c} = (1,1)^T$, $T_{end} = 1.0$ and $\varepsilon = 0$. Two temporal derivatives of the DG scheme are needed. Time step size for the smallest triangular mesh, consisting of two elements, was chosen to be $\Delta t = 0.25$. Results are computed using the multiderivative collocation method (9) with Butcher tableaux (12).

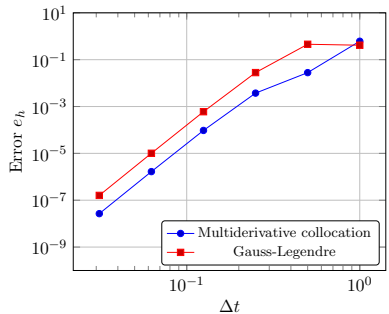

Figure 5: Numerical results for the convection equation with parameters $\mathbf{c} = (1,1)^T$, $T_{end} = 1.0$ and $\varepsilon = 0$. The polynomial order was chosen to be $p = 5$, and we compare the multiderivative collocation method against a Gauss-Legendre method. Both methods have formal order of six. Unlike in the previous figures, $\Delta t = 1.0$ was chosen on the coarsest mesh, consisting of two triangular elements.

## 5.2 Convection-diffusion equation

After having discussed the convection equation, we now turn to the convection-diffusion equation. If not stated otherwise, the interior penalty parameter is chosen to be $\eta = 20$ for a polynomial order of $p < 4$, and $\eta = 30$ for a polynomial order of $p = 5$. The numerical results being shown are all made for the equation characterized by the parameters $\mathbf{c} = (1, 1)^T$, $\varepsilon = 0.1$ and a final time of $T_{end} = 1$.

Because it is rather subtle to obtain interesting exact solutions of the convection-diffusion equation, we make use of the method of manufactured solutions and define a source term $g$ so that

$$u(x, y, t) = e^{-t} \sin(2\pi(x - t)) \sin(2\pi(y - t)) \tag{17}$$

is the exact solution to the equation. This has the added benefit of testing the ability of our algorithm (as well as our code) to handle source terms. However one drawback, from a practical point of view, is that $\partial_t g$ and $\partial_{tt} g$ also need to be computed. This can become extremely tedious, but the process can be simplified via symbolic software tools.

### 5.2.1 The convection-diffusion equation: Two-point schemes

Beginning again with the two-point schemes of order three and four, respectively, we show numerical results in Fig. 6 (third-order scheme on the left, fourth-order scheme on the right). Numerical results are computed with time step of size $\Delta t = 0.5$ on the coarsest mesh that consisted of only two triangular elements. The time step size is halved in each refinement, and the spatial mesh is uniformly refined. Our results indicate that the optimal orders are obtained. The plots look similar to those of the previous section. Note that we do not compute the $p = 0$ case, because SIPG is not meaningful for that. See also Rem. 1.

Also for this test problem, we have computed the cases for other choices of $\Delta t$, and again, we find no stability issues. This is again different to the results obtained in [30], where we in particular had problems with the stability of the diffusion terms. In this work we circumvent that issue by defining a method that is equivalent to directly differentiating the method-of-lines formulation of the PDE. That, coupled with the fact that all the solvers we consider in this work are $A$-stable, leads to a stable numerical method. This is different than most Lax-Wendroff type of discretizations, where higher derivatives are typically computed using a different method than what is performed for the first derivative.

Fig. 7 shows a comparison of the two-point two-derivative schemes against the DIRK schemes already mentioned in Sec. 5.1. The picture is the same as before: the methods behave quite similarly with sometimes a slight advantage for the two-derivative schemes.

### 5.2.2 The convection-diffusion equation: Higher-order derivatives

We continue with the higher order two-point time integration schemes of order five and six, respectively, see Tbl. 1. Numerical results are presented in Fig. 8, again optimal order is obtained and stability problems have not been observed. However, for the $p = 5$ case, we find that GMRES occasionally fails to converge. This happens also for standard DIRK schemes, and even for implicit Euler. In some sense, this is to be expected, as the stiffness matrices become increasingly stiff with higher polynomial degree. We find that the direct solver included in PETSc solves the issue, however, this needs to be fixed in the future, as a direct solver is in general not feasible. Ad-hoc solution ideas include the initialization of the GMRES routine with the outcome of a lower order method.

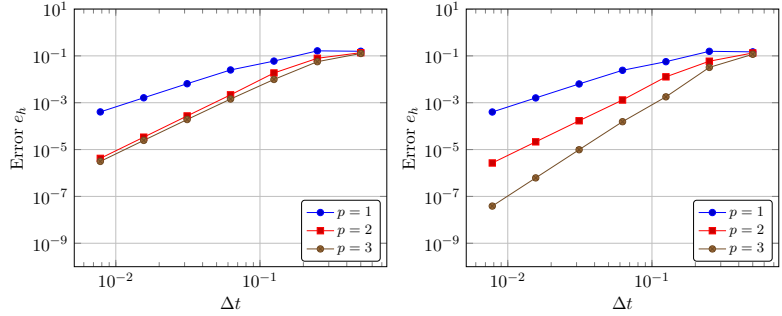

Figure 6: Numerical results for the convection-diffusion equation with parameters $\mathbf{c} = (1,1)^T$, $T_{end} = 1.0$ and $\varepsilon = 0.1$. Two temporal derivatives of the DG scheme were needed. Time step size for the smallest triangular mesh, consisting of two elements, was chosen to be $\Delta t = 0.5$. The SIPG parameter $\eta$ was chosen to be 20. Left plot: results using the third-order two-point scheme, see Tbl. 1; right plot: results using the fourth-order two-point scheme, see Tbl. 1.

### 5.2.3 The convection-diffusion equation: A three point scheme

As in the convective case, we conclude the section by showing results for the multiderivative collocation method given in (9) with Butcher tableaux as in (12). Numerical results can be seen in Fig. 9. The optimal order of $\min\{p+1, 6\}$ is attained. A comparison of the multiderivative collocation scheme against a sixth-order Gauß-Legendre scheme is made in Fig. 10. Again, we can observe a slight advantage for the multiderivative collocation scheme.

## 6 Conclusion and outlook

In this work, we introduced fully implicit multiderivative time integrators as a mechanism for discretizing convection-diffusion equations with the discontinuous Galerkin (DG) spatial discretization. Unlike most versions of DG discretizations that make use of higher derivatives such as Lax-Wendroff DG solvers, we step back and construct a solver that is equivalent to discretizing the original method-of-lines formulation of the PDE. In doing so, we sacrifice favorable properties such as being able to locally define all of our spatial operators, but the benefit of doing so includes being able to define methods that can take arbitrarily large time steps. In addition, we are able to define a single (scalar) quantity that is used to define mixed derivatives of the unknown in order to streamline the implementation and reduce the computational footprint of the solver. Future work in this direction includes considering the utility of using Gauß-Túran points for constructing higher order implicit solvers, as well as revisiting the original formulation and performing a fully discrete stability analysis in order to reduce the size of the computational stencil. In addition, we would like to implement the existing proposed solver to linear electromagnetic as well as transport dominated plasma applications such as the Vlasov-Poisson and Vlasov-Maxwell system of equations.

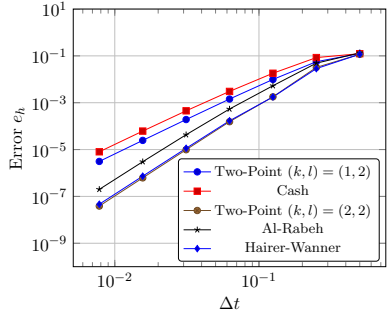

Figure 7: Numerical results for the convection-diffusion equation with parameters $\mathbf{c} = (1,1)^T$, $T_{end} = 1.0$ and $\varepsilon = 0.1$. The polynomial order was chosen to be $p = 3$, and we compare multiple time integration schemes against each other. $\Delta t = 0.5$ was chosen on the coarsest mesh, consisting of two triangular elements. The SIPG parameter $\eta$ was chosen to be 20.

## Acknowledgements

D. Seal acknowledges funding by the Naval Academy Research Council. The study of A. Jaust was supported by the Special Research Fund (BOF) of Hasselt University.

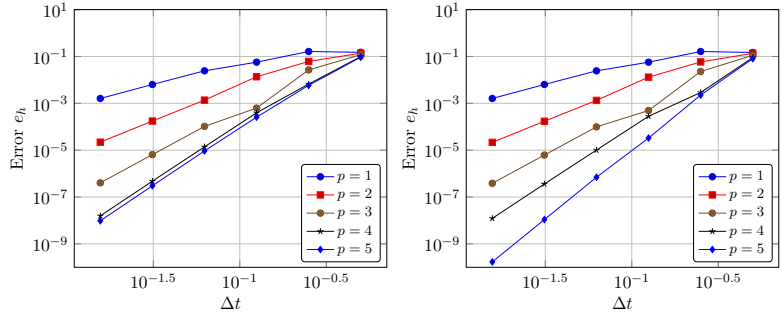

Figure 8: Numerical results for the convection-diffusion equation with parameters $\mathbf{c} = (1,1)^T$, $T_{end} = 1.0$ and $\varepsilon = 0.1$. Three temporal derivatives of the DG scheme were needed. Time step size for the smallest triangular mesh, consisting of two elements, was chosen to be $\Delta t = 0.5$. The SIPG parameter $\eta$ was chosen to be 20 ($p < 5$) and 30 ($p = 5$), respectively. Left plot: results using the fifth-order two-point scheme, see Tbl. 1; right plot: results using the sixth-order two-point scheme, see Tbl. 1.

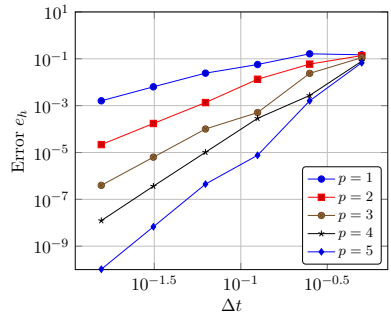

Figure 9: Numerical results for the convection-diffusion equation with parameters $\mathbf{c} = (1,1)^T$, $T_{end} = 1.0$ and $\varepsilon = 0.1$. Two temporal derivatives of the DG scheme were needed. Time step size for the smallest triangular mesh, consisting of two elements, was chosen to be $\Delta t = 0.5$. The SIPG parameter $\eta$ was chosen to be 20 ($p < 5$) and 30 ($p = 5$), respectively. Results were computed using the multiderivative collocation method (9) with Butcher tableaux (12).

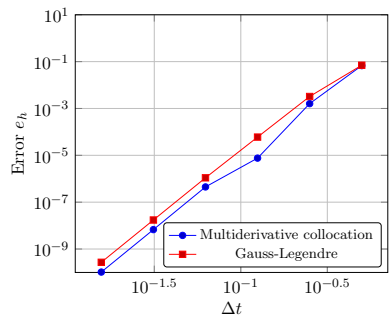

Figure 10: Numerical results for the convection-diffusion equation with parameters $\mathbf{c} = (1,1)^T$, $T_{end} = 1.0$ and $\varepsilon = 0.1$. The polynomial order was chosen to be $p = 5$, and we compare the multiderivative collocation method against a Gauss-Legendre method. Both methods have formal order of six. $\Delta t = 0.5$ was chosen on the coarsest mesh, consisting of two triangular elements. The SIPG parameter $\eta$ was chosen to be 30.